\documentclass{amsart}
\let\\=\backslash
\let\sse=\subseteq
\let\noi=\noindent
\let\wtil=\widetilde
\let\veps=\varepsilon

\def\0{\{0\}}
\def\x{{\times}}
\def\diam{{\rm diam}}
\def\size{{\rm size}}
\def\smallfrac#1#2{{\textstyle{\frac{#1}{#2}}}}
\def\Ps{{\vbox{\hbox{$\wp\kern.3pt$}\vskip.2pt}}}

\def\A{{\mathcal A}}
\def\B{{\mathcal B}}
\def\C{{\kern.5pt\mathcal C}}
\def\F{{\mathcal F}}
\def\K{{\mathcal K}}

\def\NN{{\mathbb N\kern.5pt}}
\def\RR{{\mathbb R\kern.5pt}}
\def\oRR{{\overline\RR}}

\newsymbol\triangledown 104F
\newsymbol\varnothing 203F
\newsymbol\subsetneqq 2324
\newsymbol\thicksim 2373
\let\void=\varnothing

\def\tsim{{\thicksim\,}}

\def\newmatrix#1{\null\,\vcenter{
                 \baselineskip=8pt\mathsurround=-0pt\ialign{
                 \hfil ${##}$
                 \hfil &&
                 \hfil ${##}$
                 \hfil \crcr
                 \mathstrut \crcr
                 \noalign{\kern-\baselineskip}#1 \crcr
                 \mathstrut \crcr
                 \noalign{\kern-\baselineskip} \crcr }}\!}

\newtheorem{proposition}{Proposition}

\theoremstyle{definition}

\newtheorem{remark}{Remark}

\theoremstyle{remark}

\numberwithin{theorem}{section}
\numberwithin{lemma}{section}
\numberwithin{corollary}{section}
\numberwithin{proposition}{section}
\numberwithin{conjecture}{section}
\numberwithin{definition}{section}
\numberwithin{remark}{section}
\numberwithin{question}{section}

\begin{document}

\vglue-30pt\noindent
\hfill{\it Advances in Mathematical Sciences and Applications}\/
{\bf 26} (2017) 1--18

\vglue20pt
\title[Distances Between Sets]
      {Distances Between Sets --- A Survey}
\author[A. Conci]{Aura Conci$\;\;$}
\address{Fluminense Federal University, Niteroi, Brazil}
\email{aconci@ic.uff.br}
\author[C.S. Kubrusly]{$\;\;$Carlos Kubrusly}
\address{Catholic University of Rio de Janeiro, Rio de Janeiro, Brazil}
\email{carlos@ele.puc-rio.br}
\renewcommand{\keywordsname}{Keywords}
\date{September 18, 2016}

\begin{abstract}
The purpose of this paper is to give a survey on the notions of distance
between subsets either of a metric space or of a measure space, including
definitions, a classification, and a discussion of the best-known distance
functions, which is followed by a review on applications used in many areas
of knowledge, ranging from theoretical to practical applications.
\end{abstract}

\maketitle

\section{Introduction}

Distances between subsets either of a metric space or of a measure space is
the reason of this paper, where the main focus is to define, classify, and
review the best-known distance functions, and some of their common
applications$.$ Most distances between sets, defined either for closed and
bounded (in particular, for compact, more particularly, for finite) nonempty
sets in a metric space, or for measurable sets with finite measure (in
particular, for finite sets under the counting measure) in a measure space,
become a metric themselves$.$ Behind the scenes, the main tools are naturally
the notions of metric space and measure space$.$ Although the notion of
distance between sets has a myriad of conceivable definitions (see, e.g.,
\cite[pp$.$ 46-48, 85-86, 173-184, 298-301, 359-360]{DD}, there are two
main distinct families of them, which we refer to as the
{\it Hausdorff family}\/ and the {\it Measure Theoretical family}\/$.$

\vskip6pt
Notation, terminology, basic definitions and a few results that will be
required in the sequel are posed in Section 2$.$ These are naturally bound
to the notions of metric and measure$.$ Section 3 gives a detailed account on
the Hausdorff distance (which is a metric for closed and bounded sets) and its
many relatives, which are obtained by gradual modifications of the original
\hbox{Hausdorff} distance$.$ Four equivalent forms of the original Hausdorff
measure are discussed, as well as fifteen variations of it$.$ The measure
theoretical approach, dealing with the
\hbox{Fr\'echet}--Nikod\'ym--\hbox{Aronszajn} distance (and also with its
normalized version, the Markzewisky--Steinhaus distance) is discussed in
Section 4$.$ Section 5 closes the paper with a review of the
\hbox{bibliography} (since the 80's) dealing with applications of distance
between sets towards innu\-merable subjects$.$ This is split into three
classes, namely, (1) \hbox{Computational} \hbox{Aspects} (with three
subclasses$:$ $\!${\small(1.1)} distance in graphs, {\small(1.2)} distance
between polygons, and {\small(1.3)} numerical procedures and algorithms),
(2) On Distances Between Fuzzy Sets (also with three subclasses$:$
{\small(2.1)} Markzewisky--Steinhaus distance, {\small(2.2)} \hbox{Hausdorff}
distances, {\small(2.3)} non-Hausdorff distances), and (3) Distance in Object
Analysis (again with three subsections$:$ {\small(3.1)} new metrics and
comparisons, {\small(3.2)} motion --- translation and rotation, and
{\small(3.3)} modified Hausdorff including asymmetries).

\section{Notation and Terminology}

This section summarizes classical standard topics that will be required in the
sequel, which can be found in an infinitude of books dealing with analysis in
general$.$ For instance, see \cite[Chapter 3]{EOT} for metric space properties
and \cite[Chapter 2]{MT} for measure space properties, among many others.

\vskip6pt
Let $X$ be an arbitrary nonempty set, let $d$ be a real-valued function
on the \hbox{Cartesian} product ${X\x X}$ of $X$ with itself,
$$
d\!:\!{X\x X}\!\to\RR,
$$
and consider the following properties, holding for arbitrary points
$x$, $y$, $z$ in $X$.
\begin{description}
\item{$\;$(i)}
$\;\,d(x,y)=d(y,x)\kern68pt$({\it symmetry}\/),
\item{$\,$(ii)}
$\;d(x,y)\ge0\;\;$ and $\;\;d(x,x)=0$\quad({\it nonnegativeness}\/),
\item{(iii)}
$\,d(x,y)=0\;\;$ implies $\;\;x=y\kern19pt$({\it positiveness}\/),
\item{(iv)}
$\;d(x,y)\le d(x,z)+d(z,y)\kern27pt$({\it triangle inequality}\/).
\end{description}
A real-valued function $d$ on ${X\x X}$ that satisfies all properties (i),
(ii), (iii) and (iv) is a {\it metric}\/ in $X$, and the properties themselves
are called the {\it metric axioms}\/$.$ A set $X$ equipped with a metric $d$
on ${X\x X}$ is a {\it metric space}\/, also denoted by ${(X,d)}.$ If $d$
satisfies properties (i), (ii) and (iv), but not necessarily property (iii),
then it is called a {\it pseudometric}\/ in $X.$ If $d$ satisfies properties
(i), (ii) and (iii), but not necessarily (iv), then it is sometimes called a
{\it semimetric}\/ in $X.$ What is commonly referred to as a {\it distance
function}\/ in $X$ is simply any real-valued function $d$ on ${X\x X}$ that
satisfies properties (i) and (ii) --- i.e., any symmetric nonnegative
function that vanishes at the identity line.

\begin{remark}
The difference between a metric and a pseudometric is that it is possible for
a pseudometric $d$ to vanish at a pair ${(x,y)}$ even if ${x\ne y}.$ A
pseudometric $d$ is trivial if ${d(x,y)=0}$ for every ${x,y\in X}.$ However,
given a nontrivial pseudometric $d$ in $X$ there is a natural way to obtain a
metric space ${(\wtil X,\wtil d)}$, where the set $\wtil X$ is a associated
with $X$ and $d$, and $\wtil d$ (usually denoted again by $d$) is the natural
metric in $\wtil X$ inherited from $d.$ Indeed, every nontrivial pseudometric
$d$ induces an equivalence relation $\tsim$ on $X$ (given by ${x\tsim x'}$
--- read $x$ is equivalent to $x'$ --- if and only if ${d(x,x')=0}$), and
$\widetilde X$ is the quotient space ${X/\tsim}$ which is precisely the
collection of all equivalence classes $[x]=\{x\in X\!:x'\tsim x\}$ with
respect to $\tsim$ for every element $x$ in $X$ (i.e., ${X/\tsim}$ is a
collection of sets $[x]$, called equivalence classes, such that each element
in ${X/\tsim}$ is a set consisting of all elements from $X$ that are
equivalent to each other), and the metric $\wtil d$ in $\wtil X$ is given by
$\wtil d({[x],[y]})=d({x,y})$, which does not depend on the representatives
${x\in[x]}$ and ${y\in[y]}$ from the equivalence classes.
\end{remark}

The {\it power set}\/ $\Ps(X)$ of a given set $X$ is the collection of all
subsets of $X.$ Let $X$ be an arbitrary nonempty set$.$ Equip $X$ with a
metric ${d\!:\!{X\x X}\to\RR}$ and consider the metric space ${(X,d)}.$
A nonempty subset $A$ of $X$ is {\it bounded}\/ if
$$
\sup_{x,y\,\in A}d(x,y)<\infty,
$$
otherwise $A$ is said to be unbounded, which is denoted by
$\sup_{x,y\,\in A}d({x,y})=\infty.$ The {\it diameter}\/ of a nonempty
bounded subset $A$ of $X$ is the real number
$$
\diam(A)=\!\sup_{x,y\,\in A}d(x,y).
$$
The {\it distance from a point\/ ${x\in X}$ to a nonempty set}\/
${A\in\Ps(X)}$ is the real number
$$ 
d(x,A)=\inf_{a\in A}d(x,a),
$$
and the {\it ordinary distance function between two nonempty sets $A$ and $B$
in}\/ $\Ps(X)$ is the real number
\vskip-2pt\noi
$$
d(A,B)=\!\inf_{a\in A,\,b\in B}d(a,b).
$$
The above expression defines a mere distance function
${d\!:\Ps(X)\\\void\x\Ps(X)\\\void\to\RR}.$ In fact, such a function $d$
trivially satisfies properties (i) and (ii) for all sets in ${\Ps(X)\\\void}.$
It is however clear that in general $d$ does not satisfy property (iii), and
it does not satisfy property (iv) as well --- e.g., take ${A,B,C}$ in
${\Ps(X)\\\void}$ (which may even be pairwise disjoint) such that
${d(A,B)\ne0}$ and $d({A,C})={d(C,B)=0}$.

\vskip6pt
A {\it $\sigma$-algebra}\/ $\A(X)$ of subsets of a given set $X$ is a
subcollection of the power set, ${\A(X)\sse\Ps(X)}$ (not necessarily
a proper subcollection), such that
\begin{description}
\item{$\kern12pt$}
$\;$the whole set $X$ and the empty set $\void$ belong to $\A(X)$,
\item{$\kern12pt$}
$\;$the complement ${X\\E}$ of a set $E$ in $\A(X)$ belongs to $\A(X)$,
\item{$\kern12pt$}
$\;$the union of a countable collection of sets in $\A(X)$ belongs to $\A(X)$.
\end{description}
Sets in $\A(X)$ are called {\it measurable sets}\/, and the pair
${(X,\A(X))}$ consisting of a set $X$ and a $\sigma$-algebra of subsets of
it is referred to as a {\it measurable space}\/$.$ A {\it measure}\/ is an
extended real-valued function $\mu$ on a $\sigma$-algebra $\A(X)$,
$$
\mu\!:\A(X)\to\oRR,
$$
where ${\oRR=\RR\cup\{-\infty\}\cup\{+\infty\}}$ stands for the extended real
line, satisfying the fol\-lowing properties (referred to as the
{\it measure axioms}\/).
\begin{description}
\item{(a)}
$\;{\mu(\void)=0}$,
\item{(b)}
$\;{\mu(E)\ge0}\;\;$ for every $\;\;{E\in\A(X)}$,
\item{(c)}
$\;\mu\big(\bigcup_nE_n\big)=\sum_n\mu(E_n)$
\vskip1pt\noi
\hskip6pt 
for every countable family $\{E_n\}$ of pairwise disjoint sets in $\A(X)$.
\end{description}
A {\it measure space}\/ is a triple ${(X,\A(X),\mu)}$ consisting of an
arbitrary set $X$, a $\sigma$-algebra $\A(X)$ of subsets of $X$, and a
measure $\mu$ on $\A(X).$ We assume throughout this paper that all
measures are nonzero (i.e., ${\mu(X)>0}$).

\vskip6pt
Consider an arbitrary nonempty subcollection $\C(X)$ of $\Ps(X).$ If an
arbitrary distance function $d$ on ${\C(X)\\\void\x\C(X)\\\void}$ is intended
for gauging resemblance of sets in ${\C(X)\\\void}$, then it may be useful to
control it by preventing too large (and too small) values$.$ In this case $d$
can be normalized in terms of a measure $\mu$ on a $\sigma$-algebra of subsets
of $X$, bringing forth new distance functions$.$ For instance, suppose $\mu$
is any nonzero {\it finite}\/ measure (i.e., ${0<\mu(X)<\infty}$) on any
$\sigma$-algebra $\A(X)$ of subsets of $X$, let
${d\!:\C(X)\\\void\x\C(X)\\\void\to\RR}$ be the ordinary distance function
in ${\C(X)\\\void}$, and consider the distance function
${d{\kern1pt'}\!:\C(X)\\\void\x\C(X)\\\void\to\RR}$ given by
$$
d{\kern1pt'}(A,B)=\smallfrac{1}{\mu(X)}d(A,B)
=\smallfrac{1}{\mu(X)}\inf_{a\in A,\,b\in B}d(a,b)
$$
for every sets ${A,B}$ in ${\C(X)\\\void}$, or the distance functions in the
intersection of $\C(X)$ and $\A(X)$, say
${d{\kern1pt''}\!,d{\kern1pt'''}
\!:(\C(X)\cap\A(X))\\\void\,\x\,(\C(X)\cap\A(X))\\\void\to\RR}$,
given for every sets ${A,B}$ in ${(\C(X)\cap\A(X))\\\void}$ by
$$
d{\kern1pt''}(A,B)=\smallfrac{1}{\mu(A\cup B)}d(A,B)
\quad\;\hbox{and}\;\quad
d{\kern1pt'''}(A,B)=\smallfrac{1}{\mu(A)+\mu(B)}d(A,B)
$$
if ${\mu(A)\ne0}$ or ${\mu(B)\ne0}$, otherwise $d{\kern1pt''}(A,B)$ and
$d{\kern1pt'''}(A,B)$ are defined to be zero.

\vskip6pt
Let $\#A$ denotes the cardinality of an arbitrary set ${A\in\Ps(X)}.$ Let
${\F(X)\sse\Ps(X)}$ stand for the collection of all {\it finite} subsets of
the nonempty set $X.$ Consider the collection
${\F(X)\\\void\sse\Ps(X)\\\void}$ of all finite nonempty subsets of $X.$ Take
arbitrary finite nonempty sets $A=\{a_i\}_{i=1}^m$ and $B=\{b_i\}_{i=1}^n$ in
${\F(X)\\\void}$, where ${m=\#A}$ and ${n=\#B}$ lie in $\NN$, the set of all
positive integers$.$ Now consider the restriction
$d|_{\F(X)\\\void\x\F(X)\\\void}$ of the ordinary distance function
${d\!:\!\Ps(X)\\\void\x\Ps(X)\\\void\to\RR}$ to pairs of finite nonempty sets,
${d\!:\!\F(X)\\\void\x\F(X)\\\void\to\RR}$, denoted again by the same symbol
$d.$ Thus in this case, for every ${x\in X}$ and every ${A,B\in\F(X)\\\void}$,
$$
d(x,A)\;=\inf_{a\in A}d(x,a)=\min_{a\in A}d(x,a)=\min_{1\le i\le m}d(x,a_i),
$$
$$
d(A,B)
\;=\!\inf_{a\in A,\,b\in B}d(a,b)
\;=\!\min_{a\in A,\,b\in B}d(a,b)
\;=\!\min_{1\le i\le m\;,\;1\le j\le n}d(a_i,b_j);
$$
similarly for the normalized versions (to avoid trivialities, suppose
${X\in\F(X)})$:
$$
d{\kern1pt'}({A,B})\kern-1pt=\kern-1pt\smallfrac{1}{\#X}d(A,B),
\;\;\,
d{\kern1pt''}({A,B})\kern-1pt=\kern-1pt\smallfrac{1}{\#(A\cup B)}d(A,B),
\;\;\,
d{\kern1pt'''}({A,B})\kern-1pt=\kern-1pt\smallfrac{1}{\#A+\#B}d(A,B).
$$
Still in these cases the functions
$d,d{\kern1pt'},d{\kern1pt''},d{\kern1pt'''}$ on
${\F(X)\\\void\x\F(X)\\\void}$ are distance functions but it is clear that
properties (iii) and (iv) may fail (even if ${\#X<\infty}$).

\section{The Hausdorff Family}

Let ${(X,d)}$ be a metric space$.$ Perhaps the best-known candidate for a
metric in a subset of ${\Ps(X)\\\void}$ is the \hbox{Hausdorff} function
${h\!:\Ps(X)\\\void\x\Ps(X)\\\void\to\oRR}$, which is defined for every
${A,B\in\Ps(X)\\\void}$ by \cite[p.293]{Hau1} (see also \cite[p.167]{Hau2}).
\begin{eqnarray*}
h(A,B)
&\kern-6pt=\kern-6pt&
\max\big\{\sup_{a\in A}d(a,B)\,,\,\sup_{b\in B}d(b,A)\big\}             \\
&\kern-6pt=\kern-6pt&
\max\big\{\sup_{a\in A}\inf_{b\in B^{^{\phantom I}\!\!}}d(a,b)\,,
\,\sup_{b\in B}\inf_{a\in A^{^{\phantom I}\!\!}}d(a,b)\big\}.
\end{eqnarray*}
Although the Hausdorff function satisfies properties (i) and (ii) and has the
property that $h({\{a\},\{b\}})=d({a,b})$ for every pair of singletons
${\{a\},\{b\}\in\Ps(X)\\\void}$, it is not a distance function because it is
not real-valued$.$ For instance, on ${\Ps(\RR)\\\void\x\Ps(\RR)\\\void}$ we
get ${h(\0,\RR)}=+\infty.$ Even if it were plausible to admit an
extended-real-valued distance function, the function
${h\!:\!\Ps(\RR)\\\void\x\Ps(\RR)\\\void\kern-.5pt\to\kern-1pt\oRR}\kern-.5pt$
would not be a pseudometric since property (iii) fails (e.g., ${h(A,B)=0}$ for
${A=(0,1)}$, ${B=[0,1]}$ in $\Ps(\RR)\\\void$)$.$ An equivalent formulation
for the \hbox{Hausdorff} function (also called \hbox{Blaschke} function
\cite[p.48]{DD}, or \hbox{Pompeiu}--\hbox{Hausdorff} function
\cite[Example 4.3]{RW}, \cite{BT}) is given below.

\begin{proposition}
For every ${A,B\in\Ps(X)\\\void}$,
$$
h(A,B)=\sup_{x\in X}|d(x,A)-d(x,B)|.
$$
\end{proposition}

\noi
{\it Proof}\/.
Take ${x\in X}$, ${a\in A}$, ${b\in B}$$.$ Since
$d({x,b})\le d({x,a})+d({a,b})$, we get
$d({x,B})=\inf_{b\in B}d({x,b})\le d({x,a})+\inf_{b\in B}d({a,b})
\le d({x,a})+\sup_{a\in A}\inf_{b\in B}d({a,b})$,
and hence
$d({x,B})\le\inf_{a\in A}d({x,a})+\sup_{a\in A}\inf_{b\in B}d{(a,b})
=d({x,A})+\sup_{a\in A}\inf_{b\in B}d{(a,b}).$
Symmetrically, $d({x,A})\le d({x,B})+\sup_{b\in B}\inf_{a\in A}d{(a,b}).$
Thus, for every ${x\in X}$,
$$
|d(x,A)-d(x,B)|
\le \max\big\{\sup_{a\in A}\inf_{b\in B^{^{\phantom I}\!\!}}d(a,b)\,,
\,\sup_{b\in B}\inf_{a\in A^{^{\phantom I}\!\!}}d(a,b)\big\}.
$$
In particular, for $x={a\in A}$, we get
$\sup_{a\in A}\inf_{b\in B}d({a,b})=\sup_{x=a\in A}d(x,B)=
\sup_{x=a\in A}|d(x,A)-d(x,B)|.$
Symmetrically, by setting $x={b\in B}$, we also get 
$\sup_{b\in B}\inf_{a\in A}d({a,b})=\sup_{x=b\in B}|d(x,A)-d(x,B)|.$
Therefore,
$$
\max\big\{\sup_{a\in A}\inf_{b\in B^{^{\phantom I}\!\!}}d(a,b)\,,
\,\sup_{b\in B}\inf_{a\in A^{^{\phantom I}\!\!}}d(a,b)\big\}
\le\sup_{x\in X}|d(x,A)-d(x,B)|.                                 \eqno{\qed}
$$

\vskip4pt
Although the Hausdorff function is not even a distance function when defined
on the domain ${\Ps(X)\\\void\x\Ps(X)\\\void}$, if restricted to an
appropriate subcollection it becomes a metric$.$ Let ${\B(X)\sse\Ps(X)}$
denote the collection of all {\it closed and bounded}\/ subsets of the metric
space ${(X,d)}$, take the collection ${\B(X)\\\void\sse\Ps(X)\\\void}$ of all
nonempty closed and bounded subsets of ${(X,d)}$, and consider the
restriction $h|_{\B(X)\\\void\x\B(X)\\\void}$ of
${h\!:\Ps(X)\\\void\x\Ps(X)\\\void\to\oRR}$ to ${\B(X)\\\void\x\B(X)\\\void}$,
denoted by the same symbol,
\vskip-2pt\noi
$$
{h\!:\B(X)\\\void\x\B(X)\\\void\to\RR}.
$$
This is not only real-valued but is a metric in ${\B(X)\\\void}$
\cite[Problem IX.4.8]{Dug}, \cite[Example 1.9.7]{Blu}$.$ In particular, since
${\F(X)\sse\B(X)}$, when restricted still further, now to
${\F(X)\\\void\x\F(X)\\\void}$ (i.e., when restricted to nonempty finite
sets), the Hausdorff function (again, denoted by the same symbol $h$)
$$
h\!:\F(X)\\\void\x\F(X)\\\void\to\RR,
$$
given for every $A=\{a_i\}_{i=1}^m$ and $B=\{b_i\}_{i=1}^n$ in
${\F(X)\\\void}$ by
$$
h(A,B)
=\max\big\{\!\max_{1\le i\le m}\;\min_{1\le j\le n}d(a_i,b_j)\,,
           \,\max_{1\le j\le n}\;\min_{1\le i\le m}d(a_i,b_j)\big\},
$$
is a metric in ${\F(X)\\\void}.$ For more on theoretical aspects of Hausdorff
distance see, for instance, \cite[Chapter 2]{Sen2} (also \cite{Sen1},
\cite{Sen3}), and for the topology of \hbox{Hausdorff} distances see
\cite{ALW}$.$ Indeed, the Hausdorff metric has many appealing properties,
but it shares some practical drawbacks \cite[Section 3.4]{Bad},
\cite[Section 2]{EM}$.$ For an account on theoretical and practical
aspects of some common distances between finite sets, including Hausdorff's
see, e.g., \cite{Bad, HKR, DJ, EM, Fuj, GKDS} and the references therein.

\vskip6pt
Let $\C(X)$ stand either for $\B(X)$ or $\F(X).$ Associated with each
${A\in\Ps(X)}$
set
$$
\size(A)=\left\{\newmatrix
               {\kern-20pt\#A, & \kern-27pt {\rm if}\;\; \C(X)=\F(X),    \cr
                               &                                         \cr
                \diam(A),      & \quad {\rm if}\;\; \C(X)=\B(X)\ne\F(X), \cr}
         \right.
$$
(to avoid trivialities suppose ${\size(X)>0}$) and consider the Hausdorff
metric 
$$
{h\!:\C(X)\\\void\x\C(X)\\\void\to\RR}.
$$
A natural normalization ${h'\!:\C(X)\\\void\x\C(X)\\\void\to\RR}$ of $h$ is
given by
\begin{eqnarray*}
h'(A,B)
=\smallfrac{1}{\size(X)}\,h(A,B)
&\kern-6pt=\kern-6pt&
\smallfrac{1}{\size(X)}
\max\big\{\sup_{a\in A}d(a,B)\,,\,\sup_{b\in B}d(b,A)\big\}              \\
&\kern-6pt=\kern-6pt&
\smallfrac{1}{\size(X)}\,\sup_{x\in X}|d(x,A)-d(x,B)|
\end{eqnarray*}
(cf$.$ Proposition 3.1), which is again a metric (since $h$ is a metric) if
${\size(X)<\infty}.$ When dealing with normalized versions of metrics on
${\C(X)\\\void\x\C(X)\\\void}$ it is advisable to assume that
${\size(X)<\infty}$ (i.e., $X$ is bounded if $\C(X)=\B(X)\ne\F(X)$, or
$X$ is finite if $\C(X)=\F(X)$) in order to avoid trivial pseudometrics$.$
Further natural normalizations of $h$, say
${h'',h'''\!:\C(X)\\\void\x\C(X)\\\void\to\RR}$, lead to the distances
$$
h''(A,B)=\smallfrac{1}{\size(A\cup B)}\,h(A,B)
\quad\;\hbox{and}\;\quad
h'''(A,B)=\smallfrac{1}{\size(A)+\size(B)}\,h(A,B)
$$
if ${\size(A)\ne0}$ or ${\size(B)\ne0}$; otherwise they are zero$.$
Straightforward modifications of $h$, namely,
${\wtil h,\wtil h'\!:\C(X)\\\void\x\C(X)\\\void\to\RR}$,
where maximum is replaced with sum,
$$
\wtil h (A,B)=\sup_{a\in A}d(a,B)+\sup_{b\in B}d(b,A),
$$
\vskip-2pt\noi
$$
\wtil h'(A,B)=\smallfrac{1}{\size(X)}\,\wtil h(A,B)
=\smallfrac{1}{\size(X)}\Big(\sup_{a\in A}d(a,B)+\sup_{b\in B}d(b,A)\Big),
$$
are still metrics in ${\C(X)\\\void}$ (if ${\size(X)<\infty}$
--- see \cite[Section 1]{Sho})$.$ Again, further natural normalizations
${\wtil h'',\wtil h'''\!:\C(X)\\\void\x\C(X)\\\void\to\RR}$ of $\wtil h$ are
$$
\wtil h''(A,B)=\smallfrac{1}{\size(A\cup B)}\,\wtil h(A,B)
\quad\;\hbox{and}\;\quad
\wtil h'''(A,B)=\smallfrac{1}{\size(A)+\size(B)}\,\wtil h(A,B)
$$
if ${\size(A)\ne0}$ or ${\size(B)\ne0}$, otherwise the distances are set down
to zero.

\vskip6pt
Another modified version ${h_p\!:\F(X)\\\void\x\F(X)\\\void\to\RR}$ of $h$,
where the supremum is replaced with a $p$-sum, given by
$$
h_p(A,B)=\Big({\sum}_{x\in X}|d(x,A)-d(x,B)|^p\Big)^{\frac{1}{p}}
$$
for any real ${p\ge1}$ if ${X\in\F(X)}$, and its normalization
${h_p'\!:\F(X)\\\void\x\F(X)\\\void\to\RR}$,
$$
h_p'(A,B)
=\Big(\smallfrac{1}{\#X}{\sum}_{x\in X}|d(x,A)-d(x,B)|^p\Big)^{\frac{1}{p}},
$$
are referred to as $p$-Hausdorff distances (see \cite[p.48,359]{DD}), which
in fact are metrics \cite[eq.21]{Bad})$.$ Further variants of $h$, viz.,
${h_k,h_k',\!:\F(X)\\\void\x\F(X)\\\void\to\RR}$ for ${k=0,1}$
(not all metrics), read as follows$.$ First modify $h$ by replacing supremum
with sum,
$$
h_0(A,B)= \max\Big\{{\sum}_{a\in A}d(a,B)\,,\,{\sum}_{b\in B}d(b,A)\Big\}.
$$
\vskip-4pt\noi
Its averaged version,
$$
h'_0(A,B)
=\max\Big\{\smallfrac{1}{\#A}{\sum}_{a\in A}d(a,B)\,,
\,\smallfrac{1}{\#B}{\sum}_{b\in B}d(b,A)\Big\},
$$
is called modified Hausdorff distance in \cite[p.360]{DD} (see also
\cite[eq.6,8]{DJ})$.$ Next modify $h$ by replacing
maximum with sum in $h_0$,
$$
h_1(A,B)={\sum}_{a\in A}d(a,B)+{\sum}_{b\in B}d(b,A)
$$
(considered in \cite{EM} with a multiplicative factor of $\frac{1}{2}$),
and again its averaged version,
$$
h'_1(A,B)
=\smallfrac{1}{\#A}{\sum}_{a\in A}d(a,B)
+\smallfrac{1}{\#B}{\sum}_{b\in B}d(b,A)
$$
(considered in \cite[eq.6,9]{DJ}) with a factor of $\frac{1}{2}$)$.$
Caution$:$ $h_1$ and $h_1'$ are different from $h_p$ and $h_p'$ for
${p=1}.$ A normalized version ${h_1''\!:\F(X)\\\void\x\F(X)\\\void\to\RR}$
of $h_1$, 
$$
h_1''(A,B)
=\smallfrac{1}{\#A+\#B}\Big({\sum}_{a\in A}d(a,B)+{\sum}_{b\in B}d(b,A)\Big),
$$
is referred to as geometric mean error between two images in \cite[p.360]{DD}
(see also \cite[eq.6,10]{DJ})$.$ Properties (i) and (ii) hold trivially, and
property (iii) is readily verified, and hence all the above functions are
semimetrics in ${\F(X)\\\void}.$ However, property (vi) may fail$.$ For
instance, if $A=\{a\}$, $B=\{b,c\}$, and $C=\{c\}$ with
$0<d({b,c})<$ $d({a,b})$ and $0<d({a,c})\le d({a,b})$, then
$h_k({A,C})+h_k({C,B})<h_k({A,B})$ for ${k=0,1}.$ So the semimetrics
$h_0$, $h_1$ are a not metrics in ${\F(X)\\\void}.$ On the other hand, it was
show in \cite[Theorem 1]{Fuj} that the function
${h_1'''\!:\F(X)\\\void\x\F(X)\\\void\to\RR}$ given by
$$
h_1'''(A,B)
=\smallfrac{1}{\#(A\cup B)}
\Big(\smallfrac{1}{\#A}{\sum}_{a\in A}\,{\sum}_{b\in B\\A}\!d(a,b)
+\smallfrac{1}{\#B}{\sum}_{b\in B}\,{\sum}_{a\in A\\B}\!d(a,b)\Big)
$$
is a metric in ${\F(X)\\\void}.$ For more combinations along these
lines, including theoretical aspects or practical applications, see, e.g.,
\cite{Bad, HKR, DJ, EM, DD, Fuj, GKDS}.

\vskip6pt
The {\it envelope of radius\/ $\veps$ centered at an arbitrary set}\/
${A\in\Ps(X)\\\void}$, also referred to as the $\veps$-{\it envelope}\/
(or even the $\veps$-{\it neighborhood}\/) of $A$ is the set
$$
A_\veps=\big\{x\in X\!:d(x,A)\le\veps\big\}.
$$
If $A$ is closed, so that $d(x,A)=\inf_{a\in A}d({x,a})=\min_{a\in A}d({x,a})$
for every ${x\in X}$, then $A_\veps=\{{x\in X}\!:{d(x,a)\le\veps}$ for some
${a\in A}\}.$ Another equivalent formulation for the Hausdorff function reads
as follows \cite[Problem 4.D]{Kel} (also \cite[Section 9.1]{Fal2}).

\begin{proposition}
For every ${A,B\in\Ps(X)\\\void}$,
$$
h(A,B)
=\inf\big\{\veps\ge0\!:A\sse B_\veps\;\hbox{\rm and}\;B\sse A_\veps\big\}.
$$
\end{proposition}

\noi
{\it Proof}\/.
Take ${A,B\in\Ps(X)\\\void}$ arbitrary, set
$$
h_A=\sup_{a\in A}d(a,B)
\quad\;\hbox{and}\;\quad
h_B=\sup_{b\in B}d(b,A),
$$
\vskip-4pt\noi
and observe that
$$
A\sse B_\veps \iff h_A\le\veps
\quad\;\hbox{and}\;\quad
B\sse A_\veps \iff h_B\le\veps.
$$
(Indeed, ${A\sse B_\veps}$ if and only if ${d(a,B)\le\veps}$
for all ${a\in A}).$ Thus
$$
h(A,B)=\max\big\{h_A,h_B\big\}\le
\inf\big\{\veps\ge0\!:A\sse B_\veps\;\hbox{\rm and}\;B\sse A_\veps\big\}.
$$
Conversely, note that
$$
A\sse B_{h_A}
\quad\;\hbox{and}\;\quad
B\sse A_{h_B}
$$
(in fact, ${A\sse}\{{x\in X}\!:d({x,B})\le\sup_{a\in A}d({a,B})\}= B_{h_A}).$
Therefore,
$$
\inf\big\{\veps\ge0\!:A\sse B_\veps\;\hbox{\rm and}\;B\sse A_\veps\big\}
\le\max\big\{h_A,h_B\big\}=h(A,B).                               \eqno{\qed}
$$

\vskip4pt
This version of the Hausdorff function is particularly useful to verify
that the metric space ${(\B(X)\\\void\x\B(X)\\\void,h)}$ is complete if
${X=\RR^n\!}$ equipped with its usual Euclidian metric (see, e.g.,
\cite[p.37]{Fal1}) --- this can be extended to any complete metric space
${(X,d)}$ if $\B(X)$ is swapped with $\K(X)$, the collection of all compact
subsets of the metic space ${(X,d)}$ (where
$\F(X)\sse\K(X)\sse\B(X)\sse\Ps(X)).$ Since ${A\sse B}$ if and only if
${A\\B=\void}$, for any sets $A$ and $B$, it follows that the equivalent
expression for the Hausdorff function in Proposition 3.2 can also be
rewritten as
$$
h(A,B)=
\inf\big\{\veps\ge0\!:A\\B_\veps=\void\;\hbox{\rm and}\;B\\A_\veps=\void\big\},
$$
yielding still another equivalent way to write the same \hbox{Hausdorff}
function $h$ which, if acting on ${\B(X)\\\void\x\B(X)\\\void}$ (in
particular, on ${\F(X)\\\void\x\F(X)\\\void}$), is a metric$.$ This last
form has been used in some applications (see, e.g., \cite{Con}).

\section{The Measure Theoretical Family}

Let $X$ be any nonempty set and consider the power set $\Ps(X).$ The
symmetric difference of two sets ${A,B\in\Ps(X)}$ is the set
$$
A\triangledown B=(A\\B)\cup(B\\A)=(A\cup B)\\(A\cap B),
$$
so that ${A=B}$ if and only if ${A\triangledown B=\void}$ (and
${A\triangledown\void}={\void\triangledown A}=A$ for every ${A\in\Ps(X)}).$
Let ${\mu\!:\A(X)\to\oRR}$ be an arbitrary nonzero measure (i.e.,
${\mu(X)>0}$) on an arbitrary $\sigma$-algebra $\A(X)$ of subsets of $X$ (in
particular, $\Ps(X)$ may itself be a possible $\sigma$-algebra of subsets of
$X$, depending on the measure $\mu).$ Consider the measure space
${(X,\A(X),\mu)}.$ Two sets ${E,F}$ in $\A(X)$ are {\it equivalent}\/ (or
$\mu$-{\it equivalent}\/), denoted by ${E\tsim F}$, if
${\mu(E\triangledown F)=0}.$ The relation $\tsim$ is an equivalence
relation on $\A(X).$ Define a function ${\delta\!:\A(X)\x\A(X)\to\oRR}$ by
$$
\delta(E,F)=\mu(E\triangledown F)
$$
for every ${E,F}$ in $\A(X).$ In general, the property
\begin{description}
\item{(iii)}
$\,\delta(E,F)=0\;\;$ implies $\;\;E=F$
\end{description}
fails; it fails for different sets whose symmetric difference (which is
nonempty) is of measure zero (i.e., for distinct equivalent sets)$.$ On the
other hand, the properties
\begin{description}
\item{$\;$(i)}
$\;\,\delta(E,F)=\delta(F,E)$,
\item{$\,$(ii)}
$\;\delta(E,F)\ge0\;\;$ and $\;\;\delta(E,E)=0$,
\item{(iv)}
$\;\delta(E,F)\le\delta(E,G)+\delta(G,F)$,
\end{description}
hold for every ${E,F,G\in\A(X)}$ \cite[p.138]{Nik}, \cite[Problems 3-4.9]{Hal}
(also \cite[Problem 6.5]{MT})$.$ Note that
$\delta({\void,E})=\delta({E,\void})=\mu(E)$ for every ${E\in\A(X)}.$ If the
measure $\mu$ is finite, then the function $\delta$ is real-valued$.$ Thus
assume from now on that the positive measure $\mu$ is finite (i.e.,
$0<\mu(X)<\infty).$ In this case $\delta$ is not only a distance function but
is a pseudometric in $\A(X)$.

\begin{remark}
$\!$This pseudometric ${\delta\!:\!\A(X)\x\A(X)\to\RR}$ is referred to as the
\hbox{Fr\'echet}--\hbox{Nikod\'ym}--\hbox{Aronszajn} distance
\cite[p.319]{MS}$.$ In fact, it is sound and clear that the distance is
introduced and proved to be a pseudometric in Nikod\'ym's paper
\cite[Definition 4]{Nik} as of 1930$.$ As for Aronszajn's credit, Nikod\'ym
himself \cite[p.134]{Nik} mentions a work ``to appear'' which we were not
able to trace back$.$ As for Fr\'echet's own claim of priority,
\cite[p.134]{Nik} refers to \cite{Fre1} and \cite{Fre2} as of 1921 and 1934,
respectively (see also \cite[Section 27.3.4]{Pin})$.$
\end{remark}

There are essentially two ways of normalizing this pseudometric$:$ either set
$$
\delta'(E,F)=\smallfrac{\delta(E,F)}{\mu(X)}
=\smallfrac{\mu(E\triangledown F)}{\mu(X)},
$$
or (since ${\{(E\cup F)\\(E\cap F),E\cap F\}}$ is a partition of ${E\cup F}$
so that $\mu({E\cup F})=\mu({E\triangledown F})+\mu({E\cap F)}$\,), set
$$
\delta''(E,F)=\left\{\newmatrix{
\smallfrac{\delta(E,F)}{\mu(E\cup F)}
=
\smallfrac{\mu(E\triangledown F)}{\mu(E\cup F)}
=1-\smallfrac{\mu(E\cap F)}{\mu(E\cup F)},
                                   \quad & {\rm if}\;\; \mu(E\cup F)\ne0, \cr
                                         &                                \cr
                          \kern-140pt 0, &  {\rm if}\;\; \mu(E\cup F)=0,  \cr}
              \right.
$$
for every sets ${E,F}$ in $A(X).$ Both functions are bounded by $1$ (i.e.,
${\delta'(E,F)\le1}$ and ${\delta''(E,F)\le1}$ since
${E\triangledown F}\sse{E\cup F}\sse X).$ The function
${\delta'\!:\A(X)\x\A(X)\to\RR}$ clearly is a pseudometric in $\A(X)$ (because
$\delta$ is), and it was proved in \cite[Section 1.2]{MS} that the function
${\delta''\!:\A(X)\x\A(X)\!\to\RR}$ is again a pseudometric in $\A(X)$, which
is referred to as the \hbox{Markzewisky}--\hbox{Steinhaus} distance (see also
\cite[pp.46,175,299]{DD}).

\vskip6pt
It is worth noticing that a normalized real-valued function on
${\A(X)\x\A(X)}$ given by $\frac{\mu(E\cap F)}{\mu(E\cup F)}$ if
${\mu(E\cup F)\ne0}$ is not even a distance function, since
$1=\frac{\mu(E\cap E)}{\mu(E\cup E)}\ne0$ for every $E$ in $\A(X)$ for which
${\mu(E)\ne0}$ (i.e., it is only symmetric and nonnegative but it does not
vanish at the identity line, so that property (ii) fails --- as well as
properties (iii) and (iv)).

\vskip6pt
The functions ${\delta,\delta'\!,\delta''\!:\A(X)\x\A(X)\to\RR}$ are
pseudometrics in $\A(X)\sse\Ps(X)$ (non\-trivial pseudometrics because $\mu$
is a nonzero measure), and so each of them induces a metric on the quotient
space ${\A(X)/\!\tsim}$ (of classes of equivalence of sets which differ from
each other in the same class only by a set of measure zero$:$ ${E\tsim F}$
$\kern-2.5pt\iff\kern-2pt$ $\mu({E\triangledown F})=0\,$ --- $\,$cf$.$
Remark 2.1)$.$ However, if the finite measure $\mu$ is zero only at the empty
set, then in this particular case ${E\tsim F}$ if and only if ${E=F}$, and
so ${\delta,\delta'\!,\delta''}$ are themselves metrics in the
$\sigma$-algebra ${\A(X)\sse\Ps(X)}$.

\vskip6pt
Set $\A(X)=\Ps(X)$ for a given nonempty set $X$ and consider the
counting measure ${\nu\!:\Ps(X)\to\oRR}$, which is defined by
$$
\nu(E)=\left\{\newmatrix{\#E,  & {\rm if}\;\; E\hbox{ is finite},         \cr
                               &                                          \cr
                      +\infty, & \quad {\rm if}\;\; E\hbox{ is infinite}, \cr}
       \right.
$$
for every set ${E\in\Ps(X)}.$ In this case, the only set of measure zero is
the empty set$.$ However, the counting measure is not finite if there exist
infinite subsets of $X.$ Thus suppose $X$ is finite (i.e., ${\#X<\infty}$)
so that ${\nu\!:\Ps(X)\to\RR}$ is given by
$$
\nu(E)=\#E
$$
for every subset $E$ of the finite set $X.$ Hence ${\nu(X)=\#X<\infty}$, and
$\nu$ is a {\it finite counting measure}\/ on $\Ps(X)$, and so the function
${\wtil\delta\!:\Ps(X)\x\Ps(X)\!\to\RR}$ defined by
$\wtil\delta(E,F)=\nu(E\triangledown F)$, equivalently,
$$
\wtil\delta(E,F)
=\#(E\triangledown F),
$$
for every subsets $E$ and $F$ of the finite set $X$ is indeed a metric in the
power set $\Ps(X).$ This function $\wtil\delta$ is the finite version of the
\hbox{Fr\'echet}--\hbox{Nikod\'ym}--\hbox{Aronszajn} pseudometric $\delta$
(cf$.$ Remark 4.1), now trivially a metric, whose natural normalizations
yield the metrics ${\wtil\delta',\wtil\delta''\!:\Ps(X)\x\Ps(X)\to\RR}$,
$$
\wtil\delta'(E,F)=\smallfrac{\#(E\triangledown F)}{\#X},
$$
$$
\wtil\delta''(E,F)=\left\{\newmatrix{
\smallfrac{\#(E\triangledown F)}{\#(E\cup F)}
=1-\smallfrac{\#(E\cap F)}{\#(E\cup F)},
                                    \quad & {\rm if}\;\; E\cup F\ne\void, \cr
                                          &                               \cr
                            \kern-99pt 0, & {\rm if}\;\; E\cup F=\void,   \cr}
                   \right.
$$
for every ${E,F\in\Ps(X)}.$ The function $\wtil\delta''$ is the finite
version of the \hbox{Markzewisky}--\hbox{Steinhaus} pseudometric $\delta''$
which, again, is a metric$.$ This particular metric is sometimes called
\hbox{Jaccard} distance (see, e.g., \cite[eq.4]{Fuj}, \cite[eq.5]{GKDS},
\cite[p.299]{DD} --- see also \cite[p.34]{LW}, \cite[p.174]{Gil}) after the
botanist Paul \hbox{Jaccard}, which has been referred to as \hbox{Tanimoto}
distance as well \cite[pp.46,299]{DD}, \cite[p.263]{Lip}$.$ For more on
metrics based on the cardinality of nonempty {\it finite}\/ sets see \cite{HN}.

\vskip6pt
Observe again that the normalized real-valued function on
$\,{\Ps(X)\x\Ps(X)}\,$ given by $\frac{\#(E\cap F)}{\#(E\cup F)}$ if
${E\cup F\ne\void}$ and $1$ otherwise is not a distance function$.$ In fact,
$\frac{\#(E\cap E)}{\#(E\cup E)}\ne0$ for every nonempty set $E$ in $\Ps(X)$,
and properties (ii), (iii) and (iv) fail$.$ For in\-stance, consider the sets
$E=\{e,f\}$, $\,F=\{f\}$, and $G=\{g\}$ for pairwise distinct points
${e,f,g\in X}$, so that $\frac{\#(E\cap F)}{\#(E\cup F)}\ne0$,
$\,\frac{\#(E\cap G)}{\#(E\cup G)}=\frac{\#(G\cap F)}{\#(G\cup F)}=0$, and
$\frac{\#(H\cap H)}{\#(H\cup H)}=1$ for every set $H$ in $\Ps(X).$ Such a
nondistance function has been refereed to as Jaccard similarity or Tanimoto
similarity in \cite[p.299]{DD}, and as Jaccard index in \cite[eq.4]{GKDS}.

\section{A Concise Review on Applications}

A collection of contributions on applications, classified into three
apparently distinct (but certainly not disjoint) classes, is considered in
this section$.$ Distance functions are supposed to act on appropriate
domains$.$ When we refer to a \hbox{Hausdorff} distance, it is understood
that it acts on an admissible domain that makes it well-defined (e.g., on
${\B(X)\\\void\x\B(X)\\\void}$ or, in particular, on
${\F(X)\\\void\x\F(X)\\\void}$, which make the \hbox{Hausdorff} function into
a metric)$.$ The distances and metrics of the families $h_{(s)}$ and
$\delta_{(s)}$ discussed in Sections 3 and 4 will be freely referred
throughout this section, where we will now proceed formally, omitting
theoretical details.

\vskip6pt
The forthcoming reference list bears no claim of completeness$.$ Perhaps a
complete list (if this were possible), supporting a brief review, would
become unacceptably large, leading to a dull catalog$.$ The objective
criterion for selecting the representatives in each class of the present
review was mostly based on citations; the subjective one relies on the
authors' taste.

\vskip6pt
Applications of the notion of distances between sets, in a variety of
sensible definitions, have been commonly (and naturally) used in many areas
of knowledge, ranging from theoretical to practical applications; for
instance, from computer science to biological sciences$.$ Roughly speaking,
these applications aim at shape analysis in a wide sense (and so they
encompass questions involving any sort of procedures towards sets distinction
in general)$.$ We propose a simple (perhaps too simple) and rough
classification of application areas into three classes.

\vskip9pt\noi
{\bf 5.1.}
{\bf Computational Aspects.}
All applications in this subsection deal with the \hbox{Hausdorff} family,
most of them exclusively with the plain Hausdorff distance $h$.

\vskip6pt\noi
5.1.1$.$
{\it Distance in graphs\/$.$}
Computational aspects, with emphasis in graph theory, for the Hausdorff
metric $h$ over finite sets and some of its variant distance functions as in
Section 3, were considered by Eiter and Mannila \cite{EM} in 1997$.$ These
were compared both from theoretical and computational points of view, whose
comparisons are specially tailored for applications to link distances
\cite[p.113]{EM} in graphs, some of them computed by polynomial time
algorithms$.$ In fact the Hausdorff metric is computable in polynomial time
and, in spite of its appealing properties, it is shown that it may not be
appropriate for some applications, since it does not take into account the
entire configuration of some finite sets$.$ A review on some distances in the
\hbox{Hausdorff} family is presented together with a review on the literature
up to then.

\vskip6pt\noi
5.1.2$.$
{\it Distance between polygons\/$.$}
Polygons (as well as polyhedra) are characterized by the position of a finite
number of points, in general in a finite-dimensional Euclidean space$.$ An
algorithm for computing the Hausdorff distance $h$ between a pair of convex
polygons was proposed by Atallah \cite{Ata1} in 1983$.$ This was extended by
Atallah, Ribeiro and Lifischtiz \cite{Ata2} in 1991, where algorithms for
computing some Hausdorff-type distances of two possibly overlapping and not
necessarily convex polygons was proposed$.$ Algorithms for computing Hausdorff
distance $h$ for general polyhedra represented by triangular meshes was
considered by Barton, \hbox{Hanniel}, \hbox{Elber} and Ki \cite{Bar} in 2010,
including a literature review regarding applications along this line ---
the reader is referred to the references therein$.$ For another approach,
using the Minimum Norm Duality Theorem (see, e.g. \cite[p.136]{Lue}) regarding
the ordinary distance function $d$ between convex sets in a normed space,
see \cite{Dax}.

\vskip6pt\noi
5.1.3$.$
{\it Numerical procedures and algorithms\/$.$}
The preceding subsection dealt with numerical aspects and algorithms as well,
although this may not have been the main purpose there$.$ Shonkwiler
\cite{Sho} considered in 1989 an algorithm for computing the $\wtil h$
distance between two images in linear time$.$ Huttenlocher and Kedem \cite{HK}
in 1990 and \hbox{Huttenlocher}, Kedem and Kleinberg \cite{HKK} in 1992
computed translates of the \hbox{Hausdorff} distance $h$ for subsets of the
real line and of the \hbox{Euclidian} plane, where the results were also
applied for comparing polygons under affine transformations, the main focus
being on computational speed$.$ Chew and Kedem \cite{CK} in 1998 proceeded
along the same line, considering more options for the metric $d$ on an
$n$-dimensional (${n\le3}$) real space $X$, such as the sup metric and the
$d_1$ metric in addition to the usual $d_2$ Euclidian metric$.$ Numerical
comparisons for estimating the \hbox{Hausdorff} dis\-tance $h$ between
discrete $3$-dimensional surfaces represented by triangular meshes, aiming at
the reduction of computational effort and memory usage, was considered by
\hbox{Aspert}, Santa-Cruz and Ebrahimi \cite{ASE} in 2002.

\vskip9pt\noi
{\bf 5.2.}
{\bf On Distances Between Fuzzy Sets.}
``There has been a number of papers proposing different extensions of the
Hausdorff metric to fuzzy sets$.$ None of these proposals behave as one would
intuitively expect'' \cite{Bra}.

\vskip6pt\noi
5.2.1$.$
{\it Marczewski--Steinhaus distance\/$.$}
Following the \hbox{Marczewski}--\hbox{Steinhaus} metric $\wtil\delta''$,
\hbox{Gardner}, Kanno, \hbox{Duncan}, and Selmic \cite{GKDS} proposed in 2014
an extension of $\wtil\delta''$, bounded by $\wtil\delta''$ itself, which was
shown to be a metric and suitable for applications in pattern recognition,
image processing, machine learning, and information retrieval$.$ As one would
expect, $\wtil\delta''$ would be too much of a metric to be used for fuzzy
sets, and so further measure theoretical distances are also considered$.$
Comparisons were implemented involving also the Hausdorff metric.

\vskip6pt\noi
5.2.2$.$
{\it Hausdorff distances\/$.$}
Applications along this lines involving Hausdorff-like distance functions
for fuzzy sets had been considered before by Rosenfeld \cite{Ros} in 1985
and by Chaudhuri and Rosenfeld \cite{CR1} in 1996, which were followed by
Boxer \cite{Box} in 1997, by Fan \cite{Fan} in 1998, and by Chaudhuri and
Rosenfeld \cite{CR2} in 1999$.$ A different approach where the Hausdorff
distance is used to generate further similarity measures for fuzzy sets
was considered by Hung and Yang \cite{HY} in 2004.

\vskip6pt\noi
5.2.3$.$
{\it Non-Hausdorff distances\/$.$}
A critical analysis on applications to fuzzy sets was undertaken by Brass
\cite{Bra} in 2002, where under an intriguing title ``on the nonexistence of
\hbox{Hausdorff}-like metrics for fuzzy sets'' he set about to discuss
plausible systems of metric axioms for fuzzy sets$.$ Fujita \cite{Fuj}
considered in 2013 some distance functions that can be applied to fuzzy sets,
where the Hausdorff metric $h$ and the Marczewski--Steinhaus finite-version
metric $\wtil\delta''$ are taken as starting points for yielding the original
metric $h_1'''$.

\vskip9pt\noi
{\bf 5.3.}
{\bf Distance in Object Analysis.}
By ``object analysis'' we simply mean ``image analysis'' in a very broad
sense, ranging from visible images to binary strings in general$.$ The
majority of applications of set distances focuses on problems inside this
classification, including all ranges of applications for pattern
recognition$.$ There is a very large set of references (most on
\hbox{Hausdorff} distance and its relatives) for image processing in such a
broad sense$.$ Also nonmetric distances (or nonmetric similarity functions),
meaning semimetrics (where the triangle inequality may fail), have been
considered for image analysis by \hbox{Jacobs}, \hbox{Weinshall} and
\hbox{Gdalyahu} \cite{JWG} in 2000$.$ Actually, a whole book on visual
recognition using \hbox{Hausdorff} distance by \hbox{Rucklidge} \cite{Ruc2}
has appeared in 1996, emphasizing computational aspects towards applications
on imaging processing, including concrete experiments and a large list of
references (which goes beyond set distances applied to image processing), to
which the reader is referred$.$ We comment on a shorter list (not included and
not disjoint with the above-mentioned) with a cutoff roughly after considering
some of the most cited articles (but not only) in order keep up with a
reasonable-size list.

\vskip6pt\noi
5.3.1$.$
{\it New metrics and comparisons\/$.$}
Baddeley \cite{Bad} presents in 1992 a rather detailed discussion on the
Hausdorff metric $h$ pointing out that, although theoretically attractive,
this is too sensitive a metric for image processing purposes, becoming
practically unstable$.$ Thus the metric $h_p'$ is introduced and
compared for ${p=2}$ with ``error measures of current use'', which is
essentially the metric $\delta'$ (and some asymmetric variants of it), and
comparisons involving classical synthetic images are also considered$.$ As
we have seen in 5.2.1 and 5.2.3, Fujita \cite{Fuj} in 2013, and also
\hbox{Gardner}, Kanno, Duncan, and Selmic \cite{GKDS} in 2014, introduced
original metrics as well, including discussions on the metrics $h$ and
$\wtil\delta''.$ \hbox{Dubuisson} and Jain \cite{DJ} considered in 1994 most
of the \hbox{Hausdorff} family of distances (and metrics) of Section 3
towards applications in image processing, comparing 24 combinations of them,
including comparisons involving real images$.$ They point out that $h_0'$
presents the best performance among their experiments.

\vskip6pt\noi
5.3.2$.$
{\it Motion -- translation and rotation\/$.$}
Rote \cite{Rot} proposed in (1991) an algorithm for computing the minimum
Hausdorff distance between subsets of the real line under translation$.$
This was also considered by Li, Shen and Li \cite{LSL} in 2008$.$
\hbox{Huttenlocher}, Klanderman and Rucklidge \cite{HKR}, and Huttenlocher
and Rucklidge \cite{HR}, investigated in 1993 the Hausdorff distance $h$ to
evaluate nearness between a model set and an image set, including comparisons
under translation and rigid motion (combined translation and rotation), and
an algorithm is proposed to compute these distances with examples using real
image$.$ Rucklidge \cite{Ruc1, Ruc3} in 1995 and 1997 considered a procedure
for searching space transformations of a model to match transformations that
minimize the Hausdorff distance $h$ between the transformed model and an
image, including examples$.$ Hossain, Dewan, Ahn and Chae, \cite{HDAC} also
proposed in 2012 an algorithm for computing Hausdorff-like distances with
application to moving objects$.$ Also see 5.1.3.

\vskip6pt\noi
5.3.3$.$
{\it Modified Hausdorff including asymmetries\/$.$}
Tak\'acs \cite{Tak} considered in 1998 a procedure for face matching based on
the Hausdorff family's distance $h'_0.$ This includes a penalty scheme to
ensure that images with large overlap are easily distinguished$.$ Experimental
results on a large set of face images are carried out$.$ Sim, Kwon and Park
\cite{SKP} considered in 1999 an asymmetric function associated to the
distance $h'_1$ for object matching, including simulations for comparisons
based on synthetic and real images$.$ Jesorsky, Kirchberg and Frischholz
\cite{JKF} also used in 2001 an asymmetric function associated to the distance
$h'_1$, applied for shape comparisons towards face detection, where
experiments were carried out with real images$.$ Zhao, Shi and Deng \cite{ZSD}
compared in 2005 asymmetric versions of $h$ and $h_1$ for object matching in
two-dimensional images, including experimental results.

\vskip0pt\noi
\bibliographystyle{amsplain}

\end{document}